\documentclass[a4paper,12pt,leqno]{amsart}
\usepackage{hyperref}

\begin{document}
\title[\scriptsize Entire solutions of differential-difference equations ] {Entire solutions of certain type of non-linear differential-difference equations}

\author[Li-Hao Wu,~ Ran-Ran Zhang,~ Zhi-Bo Huang]
{\scriptsize  Li-Hao Wu,~ Ran-Ran Zhang,~ Zhi-Bo Huang$^{*}$}

{\small {\small

{\small
\thanks{$^{*}$Corresponding author: Zhi-Bo Huang (huangzhibo@scnu.edu.cn)}}
 \maketitle

\vskip .5cm\par
 \maketitle {\scriptsize {\bf Abstract}. The existence of sufficiently many finite order meromorphic solutions of
a differential equation, or difference equation, or differential-difference equation, appears to be a good indicator of integrability. In this paper, we investigate  the  nonlinear differential-difference equations of form
\begin{equation*}
f(z)^{n}+L(z,f)=q(z)e^{p(z)},\eqno(*)
\end{equation*}
where $n\geq 2, L(z,f)(\not\equiv 0)$ is a linear differential-difference polynomial in $f(z)$, with small functions as its coefficients, $p(z)$ and $q(z)$ are non-vanishing polynomials. We first obtain that $n=2$ and $f(z)$ satisfies $\overline{\lambda}(f)=\sigma(f)=\deg p(z)$ under the assumption that the equation (*) possesses a transcendental entire solution of hyper order $\sigma_{2}(f)<1$. Furthermore, we give the exact form of the solutions of equation (*) when $p(z)=a, q(z)=b$, $\eta$ are constants  and $L(z,f)=g(z)f(z+\eta)+h(z)f^{'}(z)+u(z)f(z)+v(z)$  is a linear differential-difference polynomial in $f(z)$ with polynomial coefficients $g(z), h(z), u(z)$ and $v(z)$ such that $L(z,f)\not\equiv 0$ and $a b \eta\neq 0$.}
  \vskip .8cm\par

\maketitle
\numberwithin{equation}{section}
\newtheorem{theorem}{Theorem}[section]
\newtheorem{lemma}[theorem]{Lemma}
\newtheorem{proposition}[theorem]{Proposition}
\newtheorem{remark}[theorem]{Remark}
\newtheorem{definition}[theorem]{Definition}
\allowdisplaybreaks

\section{Introduction}
 \vskip .2cm\par
 Nevanlinna value distribution theory of meromorphic functions has been extensively applied to resolved growth\cite{Laine},
 value distribution\cite{Hayman, Laine}, and solvability of meromorhic solutions of linear and nonlinear differential
 equations \cite{HW, Laine, Li, LYZ}. However, meromorphic solutions of complex difference equations have been a subject of great interest in past decades, due to the application of classical Nevanlinna theory in difference by Ablowitz et. al.\cite{AHH}.
 Especially, a number of fundamental results on difference analogues of Nevanlinna value distribution  have been obtained,
 see \cite{CF1}$-$\cite{HKT}, \cite{LY1} .
  \vskip .2cm\par
 In what follows, a meromorphic function $f(z)$ is always understood to be  nonconstant and  meromorphic  in the whole complex plane $\mathbb{C}$. Concerning the value distribution of meromorphic functions, we assume that the reader is familiar with the basic
  Nevanlinna value distribution theory and its standard notations such as $m(r,f), N(r,f), T(r,f), S(r,f)$, et. al.,
  see e.g.\cite{Hayman, Laine}. In particular, for a meromorphic function $f(z)$, the notations of order
  \begin{equation*}
  \sigma(f):=\limsup\limits_{r\rightarrow \infty}\frac{\log^{+} T(r,f)}{\log r},
  \end{equation*}
  and the exponent of convergence of a-points of $f$ as
   \begin{equation*}
 \lambda(f):=\limsup\limits_{r\rightarrow \infty}\frac{\log N\left(r,\frac{1}{f-a}\right)}{\log r},
  \end{equation*}
  will appear frequently in the subsequent considerations.
   \vskip .2cm\par
   A meromorphic function
  $a(z)$ is called a small function relative to $f(z)$ if $T(r, a(z))=S(r, f)$, where
$S(r,f)$ is used to denote any quantity satisfying $S(r, f)=o(T(r,
f))$ as $r\rightarrow \infty$, possibly outside of an exceptional
set of finite logarithmic measure. Moreover, we shall use $P_{d}(f)$ to denote a differential polynomial in $f(z)$ and its derivatives $f^{'}, f^{''},\cdots$, with a total degree $d$, which has small functions relative to $f(z)$ as its coefficients. However, without confusion, we also use $P_{d}(f)$ to denote a differential-difference polynomial in $f(z)$, namely a polynomial in $f, f^{'}, f^{''},\cdots,$ and its shifts $f(z+c_{j})$, where $c_{j} (j=1,2,\cdots)$ are constants, with a total degree $d$.
 \vskip .2cm\par
 C.C.Yang\cite{Yang} considered finite order  transcendental entire solutions $f(z)$ of
 \begin{equation}\label{eq1.1}
 L(z,f)-p(z)f^{n}=h(z),
 \end{equation}
 where $L(z,f)$ denotes a linear differential polynomial in $f(z)$ with polynomial coefficients, $p(z)$ is a non-vanishing polynomial, $h(z)$ is entire and $n\geq 3.$ In particular, he showed that $f(z)$ has to be unique, unless $L(f)\equiv 0.$ After later, Heittokangas et al.\cite{HKL} investigated a slightly more general form of equation (\ref{eq1.1}), where $p(z), h(z)$ and the coefficients of $L(z, f)$ are meromorphic, and not necessarily of finite order. They showed that the method used by Yang could be modified to obtained similar uniqueness results for meromorphic solutions of this generalized equation, when $n\geq 4$. They also noted that if $n=1$ then the equation (\ref{eq1.1}) with meromorphic coefficients reduces into a linear differential equation, while if $n=2$ then (\ref{eq1.1}) contains the first and the second Painlev\'{e} differential equations and the Riccati differential equation.
  \vskip .2cm\par
Recently, several papers \cite{CY,YL,ZL} have been published regarding entire solutions of
difference and  differential-difference equations of the form
\begin{equation}\label{eq1.2}
f(z)^n + L(z,f) =h(z),
 \end{equation}
where $n\geq 2, L(z,f)$ is a linear differential-difference polynomial of $f(z)$,
and $h(z)$ is a meromorphic function of finite order. We now recall some results as follows.
  \vskip .2cm\par
  {\bf Theorem 1.A}{\cite{YL}}. Let $n\geq 4$ be an integer, $L(z,f)$ be a linear differential-difference polynomial of $f(z)$,
not vanishing identically, and $h(z)$ be a meromorphic function of finite order. Then the differential-difference equation (\ref{eq1.2})
possesses at most one transcendental entire solutions of finite order such that all coefficients of $L(z,f)$ are small
functions of $f(z)$. If such a solution $f(z)$ exists, then $f(z)$ is of the same order as $h(z)$.
\vskip .2cm\par
They also noted that if $n=3$ then the equation (\ref{eq1.2}) possesses three distinct entire solutions under certain assumptions, i.e.,
  \vskip .2cm\par
  {\bf Theorem 1.B}{\cite{YL,ZL}}. A nonlinear difference equation
\begin{equation}\label{eq1.3}
  \begin{split}
f(z)^{3}+q(z)f(z+1)=c\sin bz,
  \end{split}
\end{equation}
where $q(z)$ is a nonconstant polynomial and $b,c\in\mathbb{C}$ are nonzero constants, does not admit entire solution of
finite order. If $q(z)=q$ is a nonzero constant, then equation (\ref{eq1.3}) possesses three distinct entire solutions of
finite order, provided $b=3 \pi n$ and $q^{3}=(-1)^{n+1}\frac{27}{4}c^{2}$ for a nonzero constant $n$.
 \vskip .2cm\par
 Furthermore, they showed that if $n=2$ then the equation (\ref{eq1.2}) has no entire solution.
  \vskip .2cm\par
  {\bf Theorem 1.C}{\cite{YL}}. Let $p(z), q(z)$ be polynomials. Then a nonlinear difference equation of
\begin{equation}\label{eq1.4}
  \begin{split}
f(z)^{2}+q(z)f(z+1)=p(z)
  \end{split}
\end{equation}
has no transcendental entire solutions of finite order.
\vskip.2cm\par
Chen and Yang considered a more general form of (\ref{eq1.4}), and obtained a similar result.
\vskip.2cm\par
 {\bf Theorem 1.D}\cite{CY}. Let $p(z), h(z), g(z)$ be polynomials such that either $p(z)$ and $h(z)$ are linearly independent,
 or there is one and only one of $p(z)$ and $h(z)$ being identically equal to zero, and let $c, d_{1}, d_{2}, \lambda\in\mathbb{C}$
 be constants such that $d_{1}d_{2}\lambda \neq 0$ and $e^{\lambda c}\neq 1$. Then the differential-difference equation
 \begin{equation*}
  \begin{split}
f(z)^{2}+p(z)f(z+c)+h(z)f^{'}(z)+g(z)=d_{1}e^{\lambda z}+d_{2}e^{-\lambda z}
  \end{split}
\end{equation*}
has no entire solution of finite order.
\vskip .2cm\par
Later, X.Qi, J.Dou and L.Yang considered the nonlinear difference equation of the form
\begin{equation}\label{eq1.5}
  \begin{split}
f(z)^{n}+p(z)(\Delta_{c}f)^{m}=r(z)e^{q(z)},
  \end{split}
\end{equation}
where $\Delta_{c}f=f(z+c)-f(z)$ and $c$ is a nonzero constant, and obtained
\vskip.2cm\par
 {\bf Theorem 1.E}\cite{QDY}.Consider the nonlinear difference equation of the form (\ref{eq1.5}), where $p(z)\not\equiv 0, q(z), r(z)$ are polynomials, $n$ and $m$ are positive integers. Suppose that $f(z)$ is a transcendental entire function of finite order, not of period $c$. If $n>m$, then $f(z)$ cannot be a solution of (\ref{eq1.5}).
\vskip .2cm\par
In this paper, we consider the following nonlinear differential-difference equations of form
\begin{equation}\label{eq1.6}
  \begin{split}
f(z)^{n}+L(z,f)=q(z)e^{p(z)},
  \end{split}
\end{equation}
where $n\geq 2, L(z,f)$ is a linear differential-difference polynomial in $f(z)$, with small functions as its coefficients, $p(z)$ and $q(z)$ are non-vanishing polynomials.
\vskip .2cm\par
The reminder of this paper is organized as follows. In Section 2, we investigate the value distribution of transcendental entire solutions of equation (\ref{eq1.6}). We show that $\overline{\lambda}(f)=\sigma(f)=\deg p(z)$ if equation (\ref{eq1.6}) exactly exist a transcendental entire solution of hyper order $\sigma_{2}(f)<1$. In Section 3, we give the exact forms of transcendental entire solutions of equation (\ref{eq1.6}).

\section{Value distribution of transcendental entire solutions of differential-difference equations}
 \vskip .2cm\par
Recently, difference versions of Nevanlinna theory have been established, including the lemma of difference analogue of logarithmic derivative, difference analogue of the Clunie lemma and Mohon'ko lemma, and the second main theorem in differences, which are good tools in dealing with the value distribution of difference polynomials, and the meromorphic solutions of complex difference equations. Thus, in this section, by using difference analogues of Nenalinna theory, we investigate the value distribution of transcendental entire solutions of differential-difference equation (\ref{eq1.6}), and obtain following theorem.
  \vskip .2cm\par
  {\bf Theorem 2.1}. Let $n\geq 2$ be an integer, $L(z,f)$ be a linear differential-difference polynomial of $f(z)$,
not vanishing identically and with small functions as its coefficients,  $p(z)$ and $q(z)$ be two non-vanishing polynomials. If the differential-difference equation (\ref{eq1.6})
possesses a transcendental entire solution of hyper order $\sigma_{2}(f)<1$, then $n=2$ and $f(z)$ satisfies $\overline{\lambda}(f)=\sigma(f)=\deg p(z)$.
  \vskip .2cm\par
  {\bf Remark 2.1}. The differential-difference equation
\begin{equation*}
  \begin{split}
f(z)^{2}+h(z)f(z+\eta)-\eta e^{\eta}h(z)f'(z)+(\eta-1)e^{\eta}h(z)f(z)=z^{2}e^{2z}
  \end{split}
\end{equation*}
is solved by $f(z)=\pm z e^{z}$, but $\overline{\lambda}(f)=0, \sigma(f)=\deg p(z)=1$, where $p(z)=2 z$ and $L(z,f)=h(z)f(z+\eta)-\eta e^{\eta}h(z)f'(z)+(\eta-1)e^{\eta}h(z)f(z)\equiv 0$.
This shows that the assumption of $L(z,f)$,which is not vanishing identically in Theorem 2.1, can not be omitted.
\vskip .2cm\par
We now give some examples to show the result of Theorem 2.1 is arrived.
\vskip .2cm\par
 {\bf Example 2.1}.  The equation
 \begin{equation*}
  \begin{split}
f(z)^{2}+\frac{1}{2\pi i}z^{2}f(z+2\pi i)+\left(-\frac{1}{2\pi i}z^{2}-2z\right)f(z)=e^{2z}
  \end{split}
\end{equation*}
is solved by $f(z)=\pm z e^{z}+z$, where $p(z)=2z$ and $q(z)=1$. Obviously, $\overline{\lambda}(f)=\sigma(f)=\deg p(z)=1.$
\vskip .2cm\par
 {\bf Example 2.2}. The equation
 \begin{equation*}
  \begin{split}
f(z)^{2}+\frac{z}{e}f(z+1)+z f^{'}(z)-2f(z)-\frac{e-1}{e}z(z+1)=z^{2}e^{2z}
  \end{split}
\end{equation*}
is solved by $f(z)=\pm z e^{z}-z$, where $p(z)=2z$ and $q(z)=z^{2}$. Obviously, $\overline{\lambda}(f)=\sigma(f)=\deg p(z)=1.$
 \vskip .2cm\par
In order to prove Theorem 2.1, we need some lemmas as follows.
\vskip .2cm\par
The following Lemma 2.1 shows that non-vanishing polynomials $p(z)$ and $q(z)$  are necessary in Theorem 2.1.
    \vskip .2cm\par
  {\bf Lemma 2.1}. Under the assumption of $n$ and $L(z,f)$ in Theorem 2.1. If $p(z)$ is a constant or $q(z)\equiv 0$,
  then equation (\ref{eq1.6}) has no entire solution of hyper order $\sigma_{2}(f)<1$.
  \vskip .2cm\par
  {\bf Proof}. Contrary to our assertion, we suppose that equation (\ref{eq1.6}) has an entire solution with hyper order $\sigma_{2}(f)<1$. Since $p(z)$ is a constant or $q(z)\equiv 0$, we conclude from (\ref{eq1.6}),
 Valiron and Mohon'ko lemma\cite[Theorem 2.2.5]{Laine},  lemma of logarithmic derivative\cite[Theorem 2.3.3]{Laine}
 and its difference analogues on lemma of logarithmic
derivative\cite[Theorem 5.1]{HKT} that
 \begin{equation*}
  \begin{split}
n T(r,f)=T\left(r,q(z)e^{p(z)}-L(z,f)\right)\leq T(r,f)+S(r,f)
  \end{split}
\end{equation*}
and so
 \begin{equation*}
  \begin{split}
(n-1) T(r,f)\leq S(r,f),
  \end{split}
\end{equation*}
which contradicts our assumption that $n\geq 2$. The proof of Lemma 2.1 is approved.
   \vskip .2cm\par
  {\bf Lemma 2.2}. Let $n\geq 2$ be an integer, $L(z,f)$ be a linear differential-difference polynomial of $f(z)$,
not vanishing identically and with small functions as its coefficients,  $p(z)$ and $q(z)$ be two non-vanishing polynomials. If the differential-difference equation (\ref{eq1.6})
possesses a transcendental entire solution of hyper order $\sigma_{2}(f)<1$, then $f(z)$ satisfies $\sigma(f)=\deg p(z)$.
  \vskip .2cm\par
{\bf Proof}. Suppose that equation (\ref{eq1.6}) has an entire solution with hyper order $\sigma_{2}(f)<1$, we again conclude
from (\ref{eq1.6}), Valiron and Mohon'ko lemma\cite[Theorem 2.2.5]{Laine}, and difference analogues on lemma of logarithmic
derivative\cite[Theorem 5.1]{HKT} that
\begin{eqnarray*}
  \begin{split}
n T(r,f)&=n m(r,f)=m(r,f(z)^{n})\\
&=m\left(r, q(z)e^{p(z)}-L(z,f)\right)\\
&\leq m\left(r, q(z)e^{p(z)}\right)+m\left(r, f(z)\cdot\frac{L(z,f)-L(z,0)}{f(z)}\right)+m(r,L(z,0))\\
&\leq T\left(r, q(z)e^{p(z)}\right)+T(r,f)+S(r,f),
  \end{split}
\end{eqnarray*}
and so
\begin{equation}\label{eq2.1}
  \begin{split}
(n-1)T(r,f)\leq T\left(r, q(z)e^{p(z)}\right)+S(r,f).
  \end{split}
\end{equation}
Therefore, we get from (\ref{eq2.1}) that $\sigma(f)\leq \deg p(z)$. If $\sigma(f)< \deg p(z)$, we derive a contradiction from (\ref{eq1.6}) since
$\sigma(f^{n}(z)+q(z)f(z+1))<\deg p(z)$ and $\sigma\left(q(z)e^{p(z)}\right)=\deg p(z)$.
This yields that any transcendental entire solution of equation (\ref{eq1.6}) satisfies $\sigma(f)=\deg p(z)$. The proof of Lemma 2.2 is approved.
\vskip .2cm\par
{\bf Lemma 2.3}\cite[Theorem 2.3]{LY1}. Let $f(z)$ be a transcendental meromorphic solution of finite order $\sigma$ of a difference equation of the form
\begin{equation*}
  \begin{split}
U(z,f) P(z,f)=Q(z,f),
  \end{split}
\end{equation*}
where $U(z,f), P(z,f)$ and $Q(z,f)$ are difference polynomials such that the total degree $\deg U(z,f)=n$ in $f$  and its shifts, and let $\deg Q(z,f)\leq n$. Moreover, we assume that $U(z,f)$  contains just one term of maximal total degree in $f(z)$ and its shifts. Then for each $\varepsilon>0$,
\begin{equation*}
  \begin{split}
m(r,P(z,f))=O(r^{\sigma-1+\varepsilon})+S(r,f),
  \end{split}
\end{equation*}
possibly outside of an  exceptional set with finite logarithmic measure.
\vskip .2cm\par
{\bf Remark 2.2}. By using similar method\cite[Lemma 2.4.2]{Laine}, we note that Lemma 2.3 is still valid
if $f(z)$ is a transcendental meromorphic function with $\sigma(f)<\infty$
, $P(z,f)$ and $Q(z,f)$ are differantial-difference polynimials in $f(z)$. Moreover,
  \begin{equation*}
  \begin{split}
m(r,P(z,f))=O(\log r),
  \end{split}
\end{equation*}
possibly outside of an  exceptional set with finite logarithmic measure, if $f(z)$ is a transcendental entire function with $\sigma(f)=1$
, $P(z,f)$ and $Q(z,f)$ are differantial-difference polynimials in $f(z)$, with polynomial coefficients.
\vskip .2cm\par
Lemma 2.2 shows that any transcendental entire solution of equation (\ref{eq1.6}) must be of finite order. Furthermore, we will obtain the following lemma.
\vskip .2cm\par
  {\bf Lemma 2.4}. Let $n\geq 3$ be an integer, $L(z,f)$ be a linear differential-difference polynomial of $f(z)$,
not vanishing identically and with small functions as its coefficients,  $p(z)$ and $q(z)$ be two non-vanishing polynomials. Then  equation (\ref{eq1.6})
does not possess any transcendental entire solutions of finite order.
\vskip .2cm\par
  {\bf Proof}. Contrary to our assertion, we suppose that equation (\ref{eq1.6}) possesses a transcendental entire solution of finite order.
  \vskip .2cm\par
   Differentiating both sides of (\ref{eq1.6}), we obtain
\begin{equation}\label{eq2.2}
  \begin{split}
n f(z)^{n-1}f^{'}(z)+L^{'}(z, f)=[q^{'}(z)+q(z)p^{'}(z)]e^{p(z)}.
  \end{split}
\end{equation}
\vskip .2cm\par
By eliminating $e^{p(z)}$  from (\ref{eq1.6}) and (\ref{eq2.2}), we conclude that
\begin{equation}\label{eq2.3}
  \begin{split}
f(z)^{n-1}P(z,f)=Q(z,f),
  \end{split}
\end{equation}
where
\begin{equation*}
  \begin{split}
P(z,f)=n f^{'}(z)-\left(p^{'}(z)+\frac{q^{'}(z)}{q(z)}\right)f(z),
  \end{split}
\end{equation*}
and
\begin{equation*}
  \begin{split}
Q(z,f)=\left(p^{'}(z)+\frac{q^{'}(z)}{q(z)}\right) L(z, f)-L^{'}(z, f).
  \end{split}
\end{equation*}
\vskip .2cm\par
Now, we consider the following two cases.
\vskip .2cm\par
{\bf Case 1.} $Q(z,f)\equiv 0.$ Then we have from (\ref{eq2.3}) that
\begin{equation*}
  \begin{split}
P(z,f)=n f^{'}(z)-\left(p^{'}(z)+\frac{q^{'}(z)}{q(z)}\right)f(z)\equiv 0.
  \end{split}
\end{equation*}
This shows that $f(z)^{n}=\mu q(z)e^{p(z)}$, and so $f(z)=r(z)e^{\frac{p(z)}{n}}$, where $\mu$ is a nonzero constant
and $r(z)^{n}=\mu q(z)$. Substituting $f(z)$ into (\ref{eq1.6}), we obtain
\begin{equation}\label{eq2.4}
  \begin{split}
(\mu-1)q(z)e^{p(z)}+L(z, f)=0,
  \end{split}
\end{equation}
and so $L(z,f)\equiv 0$ if $\mu \equiv 1,$ a contradiction. If $\mu \neq 1,$ we apply  Valiron and Mohon'ko lemma to
(\ref{eq2.4}) to obtain that
\begin{equation*}
  \begin{split}
T\left(r, e^{\frac{p(z)}{n}}\right)+S(r,f)=T(r,L(z,f))=n T\left(r, e^{\frac{p(z)}{n}}\right)+S(r,f),
  \end{split}
\end{equation*}
and  again get a contradiction since $n\geq 2$.
\vskip .2cm\par
{\bf Case 2.} $Q(z,f)\not\equiv 0.$   We note that $f(z)$ is a finite order entire solution, $p(z)$ and $q(z)$ are polynomials. If $n\geq 3, $  we deduce from (\ref{eq2.3}) and Lemma 2.3 that
\begin{eqnarray*}
  \begin{split}
  &T(r,P(z,f))=m(r,P(z,f))+N(r,P(z,f))=S(r,f),\\
  &T(r,f P(z,f))=m(r,f P(z,f))+N(r,fP(z,f))=S(r,f),
  \end{split}
\end{eqnarray*}
 and so
\begin{equation*}
  \begin{split}
T(r,f)=T(r, fP(z,f)/P(z,f))\leq T(r,fP(z,f))+T(r,1/P(z,f))=S(r,f),
  \end{split}
\end{equation*}
a contradiction. Lemma 2.4 is approved.
 \vskip .2cm\par

  {\bf Proof of Theorem 2.1}. According to Lemma 2.2 and Lemma 2.4, we just need to prove $\overline{\lambda}(f)=\sigma(f)$ when $n=2$ in equation (\ref{eq1.6}).
\vskip .2cm\par
Similar to the proof of Lemma 2.4, we can rewrite (\ref{eq2.3}) as

\begin{equation}\label{eq2.5}
  \begin{split}
f(z)P(z,f)=Q(z,f),
  \end{split}
\end{equation}
where
\begin{equation}\label{eq2.6}
  \begin{split}
P(z,f)=2 f^{'}(z)-\left(p^{'}(z)+\frac{q^{'}(z)}{q(z)}\right)f(z),
  \end{split}
\end{equation}
and
\begin{equation*}
  \begin{split}
Q(z,f)=\left(p^{'}(z)+\frac{q^{'}(z)}{q(z)}\right) L(z, f)-L^{'}(z, f).
  \end{split}
\end{equation*}
\vskip .2cm\par
By using the same method in Lemma 2.4, we deduce a contradiction again when $Q(z,f)\equiv 0$. Thus, we just prove the case that $Q(z,f)\not\equiv 0$, which shows that  $L(z,f)\not\equiv 0$.
\vskip .2cm\par
We now deduce from (\ref{eq2.5}), Lemma 2.3 and Remark 2.2 that
\begin{equation}\label{eq2.7}
  \begin{split}
T(r, P(z,f))=m(r,P(z,f))+N(r, P(z,f))=S(r,f),
  \end{split}
\end{equation}
possibly outside of an  exceptional set with finite logarithmic measure.
\vskip .2cm\par
Furthermore, we conclude from (\ref{eq2.6}) and lemma of logarithmic derivative that
\begin{equation}\label{eq2.8}
  \begin{split}
m\left(r, \frac{1}{f(z)}\right)
&\leq m\left(r, \frac{1}{P(z,f)}\right) +m\left(r, 2 \frac{f^{'}(z)}{f(z)}-\left(p^{'}(z)+\frac{q^{'}(z)}{q(z)}\right)\right)\\
&\leq  m\left(r, \frac{1}{P(z,f)}\right) +S(r,f),
  \end{split}
\end{equation}
possibly outside of an  exceptional set with finite logarithmic measure.
\vskip .2cm\par
We now assert that $f(z)$ has infinitely many zeros. Otherwise, we can deduce $T(r,f)=S(r,f)$ from (\ref{eq2.7}), (\ref{eq2.8}) and the first main theory, a contradiction.
\vskip .2cm\par
Since $p(z)$ and $q(z)$ are polynomials, there are only finite common zeros between $f(z)$,  $p(z)$ and $q(z)$. Suppose that
$z_{0}$ is a zero of $f(z)$ with order $k$ such that  $p(z_{0})\neq 0$ and $q(z_{0})\neq 0$, then by (\ref{eq2.6}), $z_{0}$ is also a zero of $P(z, f)$ with order $k-1$. Thus, we have
\begin{equation}\label{eq2.9}
  \begin{split}
N\left(r, \frac{1}{f(z)}\right)\leq \overline{N}\left(r, \frac{1}{f(z)}\right) +N\left(r, \frac{1}{P(z,f)}\right)+O(\log r).
  \end{split}
\end{equation}
\vskip .2cm\par
We then yield from (\ref{eq2.7})$-$(\ref{eq2.9}) and the first main theory that
\begin{equation*}
  \begin{split}
T(r,f)
&=T\left(r, \frac{1}{f(z)}\right)+O(\log r)\\
&\leq \overline{N}\left(r, \frac{1}{f(z)}\right) +T\left(r, \frac{1}{P(z,f)}\right)+S(r,f)\\
&\leq \overline{N}\left(r, \frac{1}{f(z)}\right) +S(r,f),
  \end{split}
\end{equation*}
possibly outside of an  exceptional set with finite logarithmic measure. This yields $\sigma(f)\leq \overline{\lambda}(f)$. Therefore, we have  $\overline{\lambda}(f)=\sigma(f)$. The proof of Theorem 2.1 is approved.
 \vskip .5cm\par
 \section{Forms of transcendental entire solutions of two order differential-difference equations}
 \vskip .2cm\par
 The existence of sufficiently many finite order meromorphic solutions of a difference equation appears to be a good indicator of integrability. In this section, we present the exact forms of transcendental entire solutions of a certain type of second order differential-difference equations, and have the following result.
  \vskip .2cm\par
  {\bf Theorem 3.1}. Let $L(z,f)=g(z)f(z+\eta)+h(z)f^{'}(z)+u(z)f(z)+v(z)$ be a linear differential-difference polynomial in $f(z)$ with polynomial coefficients $g(z), h(z), u(z)$ and $v(z)$ such that $L(z,f)\not\equiv 0$, and let $a, b, \eta$ be constants such that $a b \eta\neq 0$. Then any finite order entire solution of
  \begin{equation}\label{eq3.1}
  \begin{split}
 f(z)^{2}+L(z,f)=b e^{a z}
  \end{split}
\end{equation}
must be form of
 \begin{equation}\label{eq3.2}
  \begin{split}
 f(z)=c e^{\frac{a}{2} z}+f_{0}(z),
  \end{split}
\end{equation}
where $c^{2}=b$ and $f_{0}(z)=-\frac{1}{2} \left(e^{\frac{a \eta}{2}} g(z)+\frac{a}{2}h(z)+u(z)\right)$ is a non-vanishing polynomial.
\vskip .2cm\par
The following example is listed to show that Theorem 3.1 is valid.
\vskip .2cm\par
{\bf Example 3.1}. Let $\eta$ be a nonzero constant such that $L(z,f)\not\equiv 0$ in Theorem 3.1 when
\begin{eqnarray*}
 {\begin{array}{llll}
g(z)=2 e^{-\eta}z, &h(z)=e^{-\eta}, &u(z)=-e^{-\eta}\\
a=2, &b=1, &v(z)=(2-e^{\eta})e^{-\eta} z^{2}+(2 \eta -1)e^{-\eta} z +e^{-\eta}.
\end{array}}
\end{eqnarray*}
Then  $f_{0}(z)=-z$ and the equation
\begin{equation*}
  \begin{split}
& f(z)^{2}+ 2 e^{-\eta}z f(z+\eta)+e^{-\eta}f^{'}(z)-e^{-\eta} f(z)\\
&+\left[(2-e^{\eta})e^{-\eta} z^{2}+(2 \eta -1)e^{-\eta} z +e^{-\eta}\right]=e^{2 z}
  \end{split}
\end{equation*}
has entire solutions $f(z)=\pm e^{z}-z$, which are the forms of (\ref{eq3.2}).

\vskip .2cm\par
We first give some lemmas.
\vskip .2cm\par
{\bf Lemma 3.1}\cite[Theorem 1.51]{YY} . Suppose that $n\geq 2$ and let $f_{j}(z),j=1,2,\cdots,n$ be meromorphic functions and $g_{j}, j=1,2,\cdots,n$ be entire functions such that\par
(i) $\sum\limits_{j=1}^{n} f_{j}(z)e^{g_{j}(z)}\equiv 0;$\par
(2) when $1\leq j<k\leq n, g_{j}(z)-g_{k}(z)$ is not constant;\par
(3) when $1\leq j\leq n, 1\leq h<k\leq n,$
\begin{equation*}
  \begin{split}
T(r,f_{j}(z))=o\{T(r, \exp\{g_{h}(z)-g_{k}(z)\})\}\quad r\rightarrow \infty, r\not\in E,
  \end{split}
\end{equation*}
where $E\subset (1,\infty)$ is of finite linear measure or finite logarithmic measure.
Then $f_{j}(z)\equiv 0, j=1,2,\cdots,n.$
\vskip .2cm\par
{\bf Lemma 3.2}. Let $a$ be a nonzero constant, and $H(z)$ be a non-vanishing  polynomial. Then the differential equation
\begin{equation}\label{eq3.3}
  \begin{split}
2 f^{'}(z)-a f(z)=H(z)
  \end{split}
\end{equation}
has a special solution $f_{0}(z)$ which is a non-vanishing polynomial.
\vskip .2cm\par
{\bf Proof}. If $H(z)$ is a nonzero constant, then clearly $f_{0}(z)=-\frac{H(z)}{a}$ is a special solution of (\ref{eq3.3}). Thus, we now suppose that
\begin{equation*}
  \begin{split}
H(z)= \lambda_{n} z^{n} + \lambda_{n-1}z^{n-1} +\cdots + \lambda_{1} z + \lambda_{0},
  \end{split}
\end{equation*}
where $n\geq 1$ is an integer, and $\lambda_{n}(\neq 0), \lambda_{n-1}, \cdots, \lambda_{0}$ are constants.
\vskip .2cm\par
We use the method of undetermined coefficients, to derive the polynomial solution
$f_0(z)$ satisfying (3.3) by $a, \lambda_{n} , \lambda_{n-1},  \cdots, \lambda_{0}$.
Clearly,  we see from (\ref{eq3.3}) that $\deg f_{0} =\deg H$.
If $n=1$,  equation (\ref{eq3.3}) has a polynomial solution
\begin{equation*}
  \begin{split}
f_0(z)=-\frac{\lambda_{1}}{a}z+\left(-\frac{\lambda_{0}}{a}-2\frac{\lambda_{1}}{a^{2}}\right).
  \end{split}
\end{equation*}
\vskip .2cm\par
If $n\geq 2$, a general case, equation (\ref{eq3.3}) has a polynomial solution
\begin{equation*}
  \begin{split}
f_0(z)= b_{n} z^{n} + b_{n-1}z^{n-1} +\cdots + b_{1} z + b_{0},
  \end{split}
\end{equation*}
where
\begin{equation*}
  \begin{split}
b_n=-\frac{\lambda_n}{a},~~b_{j}=\frac{2(j+1)b_{j+1}-\lambda_{j}}{a},~ j=0, 1,\cdots, n-1.
  \end{split}
\end{equation*}
\vskip .2cm\par
Thus, equation (\ref{eq3.3}) has a non-vanishing polynomial solution $f_0(z)$.

\vskip .2cm\par
{\bf Proof of Theorem 3.1}.  Suppose that $f(z)$ is an entire solution of equation (\ref{eq3.1}) with finite order. Similar to the proof of Lemma 2.4, we can obtain
\begin{equation}\label{eq3.4}
  \begin{split}
f(z)P(z,f)=Q(z,f),
  \end{split}
\end{equation}
where
\begin{equation*}
  \begin{split}
P(z,f)=2 f^{'}(z)-a f(z),~~~~Q(z,f)=a L(z, f)-L^{'}(z, f).
  \end{split}
\end{equation*}
\vskip .2cm\par
We now discuss the following two cases.
\vskip .2cm\par
{\bf Case 3.1}. $Q(z,f)\equiv 0$. Then equation (\ref{eq3.4}) implies that $P(z,f)=2 f^{'}(z)-a f(z)\equiv 0$, which yields
\begin{equation}\label{eq3.5}
  \begin{split}
f(z)=c e^{\frac{a}{2}z}
  \end{split}
\end{equation}
for some non-zero constant $c$. We now substitute (\ref{eq3.5}) into (\ref{eq3.1}), and conclude that
\begin{equation}\label{eq3.6}
  \begin{split}
(c^{2}-b) e^{a z}+ c\left(e^{\frac{a \eta}{2}} g(z)+\frac{a}{2}h(z)+u(z)\right)e^{\frac{a}{2}z}+v(z)=0.
  \end{split}
\end{equation}
Thus, we deduce from Lemma 3.1 and (\ref{eq3.6}) that
\begin{equation*}
  \begin{split}
  c^{2}=b,~~e^{\frac{a \eta}{2}} g(z)+\frac{a}{2}h(z)+u(z)\equiv 0~~and~~v(z)\equiv 0,
  \end{split}
\end{equation*}
which yield $L(z,f)\equiv 0$, a contradiction.
\vskip .2cm\par
{\bf Case 3.2}. $Q(z,f)\not\equiv 0$. We first obtain from Theorem 2.1 that $\sigma(f)=1$. Thus, we further apply Lemma 2.3 and Remark 2.2 to (\ref{eq3.4}) that
\begin{equation}\label{eq3.7}
  \begin{split}
m(r, 2 f^{'}(z)-a f(z))=m(r,P(z,f))=O(\log r),
  \end{split}
\end{equation}
possibly outside of an exceptional set with finite logarithmic measure.  (\ref{eq3.7}) implies that $2 f^{'}(z)-a f(z)$ is a polynomial. Therefore, we have from  (\ref{eq3.4}) and $Q(z,f)\not\equiv 0$ that
\begin{equation}\label{eq3.8}
  \begin{split}
2 f^{'}(z)-a f(z)=H(z),
  \end{split}
\end{equation}
where $H(z)$ is a nonvanishing polynomial. Thus, we obtain from Lemma 3.2 that the equation (\ref{eq3.8}) must have a non-vanishing polynomial solution, say, $f_{0}(z)$.
\vskip .2cm\par
Since the differential equation
\begin{equation*}
  \begin{split}
2 f^{'}(z)-a f(z)=0,
  \end{split}
\end{equation*}
has a fundamental solution $f(z)=e^{\frac{a}{2} z}.$ It follows that the general solution $f(z)$ of (\ref{eq3.8}) can be express as
\begin{equation}\label{eq3.9}
  \begin{split}
f(z)=c e^{\frac{a}{2} z}+f_{0}(z),
  \end{split}
\end{equation}
where $c$ is a nonzero constant, $f_{0}(z)$ is a special non-vanishing polynomial solution.
\vskip .2cm\par
Substituting (\ref{eq3.9}) into (\ref{eq3.1}) , we conclude
\begin{equation}\label{eq3.10}
  \begin{split}
&(c^{2}-b)e^{a z}+c\left(2 f_{0}(z)+ e^{\frac{a\eta}{2}}g(z)+\frac{a}{2}h(z)+u(z)\right)e^{\frac{a}{2} z}\\
&+f_{0}^{2}(z)+g(z) f_{0}(z+\eta)+h(z)f_{0}^{'}(z)+u(z)f_{0}(z)+v(z)=0.
  \end{split}
\end{equation}
\vskip .2cm\par
It follows from Lemma 3.1 and (\ref{eq3.10}) that
\begin{equation*}
  \begin{split}
  &c^{2}=b,\\
  &2 f_{0}(z)+ e^{\frac{a\eta}{2}}g(z)+\frac{a}{2}h(z)+u(z)\equiv 0,\\
  &f_{0}^{2}(z)+g(z) f_{0}(z+\eta)+h(z)f_{0}^{'}(z)+u(z)f_{0}(z)+v(z)\equiv0.
  \end{split}
\end{equation*}
\vskip .2cm\par
We further conclude that
\begin{equation*}
  \begin{split}
  f_{0}(z)=-\frac{1}{2}\left( e^{\frac{a\eta}{2}}g(z)+\frac{a}{2}h(z)+u(z)\right),
  \end{split}
\end{equation*}
and, in this case,
\begin{equation*}
  \begin{split}
L(z,f)=-  f_{0}(z)\left(  f_{0}(z)+2 c e^{\frac{a}{2}z}\right)\neq 0,
  \end{split}
\end{equation*}
since $a\neq 0$ and $  f_{0}(z)$ is a non-vanishing polynomial.
\vskip .2cm\par
Thus, any finite order entire solution of the equation (\ref{eq3.1})
must be form of (\ref{eq3.2}). The proof of Theorem 3.1 is approved.

\vskip .5cm\par

{\bf Acknowledgements:}\par The first author is supported in part by Foundation for Young Talents in General Higher
Education of Guangdong(No:2015KQNCX230).
The Second author is supported in part by Guangdong National Natural Science Foundation(No.2016A030313745).
 The third author is supported in part by Guangdong National Natural Science Foundation(No.2018A030313508).

\vskip .2cm\par

\vskip.5cm\par\quad\\
Li-Hao Wu\\
School of Computer Engineering\\
Gungzhou College of South China University of Technology\\
Guangzhou 510800, P.R.China\\
Email: wulh@gcu.edu.cn
\vskip.5cm\par\quad\\
Ran-Ran Zhang\\
Department of Mathematics\\
Guangdong University of Education\\
Guangzhou, 510303, P.R.China\\
Email:  zhangranran@gdei.edu.cn
\vskip.5cm\par\quad\\
Zhi-Bo Huang\\
School of Mathematical Sciences\\
South China Normal University\\
Guangzhou, 510631, P.R.China\\
Email:  huangzhibo@scnu.edu.cn,~zhibo.huang@ucl.ac.uk

\begin{thebibliography}{99}
\scriptsize

\bibitem {AHH} M.J.Ablowitz, R.Halburd and B.Herbst, \emph{On
the extension of Painlev\'{e} property to difference equations},
Nonlinearity 13(2000), no.3, 889-905.

\bibitem{CY} Z.X.Chen and C.C.Yang, \emph{On entire solutions of certain type of
nonlinear differential-difference equations}, Taiwanese J.Math.,18(2014), no.3,677-685.


\bibitem{CF1} Y.M.Chiang and S.J.Feng, \emph{On the Nevanlinna
characteristic of $f(z+\eta)$ and difference equations in the
complex plane}, Ramanujan J. 16(2008), no. 1, 105-129.

\bibitem{HK1} R.G.Halburd and R.J.Korhonen, \emph{Difference analogue
of the lemma on the logarithmic derivative with applications to
difference equations}, J. Math. Anal. Appl. 314(2006), no. 2,
477-487.

\bibitem{HK2} R.G.Halburd and R.J.Korhonen, \emph{Nevanlinna theory for
the difference operator}, Ann. Acad. Sci. Fenn. Math. 31(2006), no.
2, 463-478.


\bibitem{HKT} R.G.Halburd, R.J.Korhonen and K.Tohge, \emph{Holomorphic curves with shift invariant hyper plane preimages}, Trans. Amer. Math. Soc., 366(2014), no.8, 4267-4298.

\bibitem {HW} R.Halburd and J.Wang,\emph{All admissible meromorphic solutions of Hayman's equation}, Internatonal Mathematics Research Notices, 2015(2015), no.18, 8890-8902.

\bibitem {Hayman} W.K.Hayman, \emph{Meromorphic Functions}, Clarendon Press, Oxford, 1964.

\bibitem{HKL}  J.Heittokangas, R.Korhonen, I.Laine,  \emph{On meromorphic solutions of  certain  nonlinear differential equations}, Bull. Austral. Math. Soc., 66(2002), 331-343.

\bibitem {Laine} I.Laine,\emph{Nevanlinna theory and complex differental equations}, de Gruyter, Berlin, 1993.

\bibitem{LY1} I.Laine and C.C.Yang, \emph{Clunie  theorems for
difference and q-difference polynomials}, J. London Math. Soc. 76
(2007), no. 3, 556-566.

\bibitem{Li} P.Li, \emph{Entire solutions of certain type of differential equations}, J.Math.Anal.Appl., 344(2008),no.1, 253-259.

\bibitem{LYZ} L.W.Liao, C.C.Yang and J.J.Zhang, \emph{On meromorphic solutions of certain type of nonlinear
differential equations}, Ann. Acad. Sci. Fenn.Math., 38(2013),no.2, 581-593.

\bibitem{QDY} X.Qi, J. Dou, L.Yang, \emph{The properties of solutions of certain type of
difference equations}, Adv. Difference Equ.,(256)2014.

\bibitem{Yang} C.C.Yang, \emph{On entire solutions of a certain type of nonlinear differential equation}, Bull. Austral. Math. Soc., 64(2001), 377-380.

\bibitem{YL}C.C.Yang and I. Laine, \emph{On analogies between nonlinear
difference and differential equations}, Proc. Japan Acad. Ser. A,
Math. Sci. 86(2010), no. 1, 10-14.

\bibitem{YY} C.C.Yang and H.X.Yi, \emph{Uniqueness Theory of
Meromorphic Function}, Kluwer, Dordrecht, 2003.


\bibitem{ZL} J.Zhang and L.W.Liao, \emph{On entire solutions of a certain type of nonlinear differential
and difference equations}, Taiwanese J.Math.,15(2011), no.5,2145-2157.

\end{thebibliography}
\end{document}